%% file: main.tex
\documentclass{article}

\usepackage{resources/PRIMEarxiv}

\usepackage[utf8]{inputenc} % allow utf-8 input
\usepackage[T1]{fontenc}    % use 8-bit T1 fonts
\usepackage{hyperref}       % hyperlinks
\usepackage{url}            % simple URL typesetting
\usepackage{booktabs}       % professional-quality tables
\usepackage{amsfonts}       % blackboard math symbols
\usepackage{nicefrac}       % compact symbols for 1/2, etc.
\usepackage{microtype}      % microtypography
\usepackage{lipsum}
\usepackage{soul,xcolor}
\usepackage{blkarray}

\usepackage{caption}
\usepackage{subcaption}
\usepackage{float}
\usepackage{amsmath}
\usepackage{mathtools}
\usepackage{fancyhdr}       % header
\usepackage{graphicx}       % graphics
\graphicspath{{media/}}     % organize your images and other figures under media/ folder

\title{SMS: Spiking Marching Scheme for Efficient Long Time Integration
of Differential Equations}

\author{
  Qian Zhang\thanks{The first and second author contributed equally.} , Adar Kahana$^*$, George Em Karniadakis \\
  Division of Applied Mathematics \\
  Brown University \\
  Providence, RI \\
  \texttt{\{qian\_zhang1, adar\_kahana,  george\_karniadakis\}@brown.edu}  \\
   \And
  Panos Stinis \\
  Advanced Computing, Mathematics and Data Division \\
  Pacific Northwest National Laboratory \\
  Richland, WA \\
  \texttt{panagiotis.stinis@pnnl.gov} \\
}

\begin{document}
\maketitle

\begin{abstract}
We propose a Spiking Neural Network (SNN)-based explicit numerical scheme for long time integration of time-dependent Ordinary and Partial Differential Equations (ODEs, PDEs). The core element of the method is a SNN, trained to use spike-encoded information about the solution at previous timesteps to predict spike-encoded information at the next timestep. After the network has been trained, it operates as an explicit numerical scheme that can be used to compute the solution at future timesteps, given a spike-encoded initial condition. A decoder is used to transform the evolved spiking-encoded solution back to function values. We present results from numerical experiments of using the proposed method for ODEs and PDEs of varying complexity.
\end{abstract}


\input{introduction}
\input{background}
\input{sms}
\input{results}
\input{conclusion}
\input{acknowledgments}

\bibliographystyle{unsrt}  
\bibliography{references}

\end{document}

%% file: introduction.tex
\section{Introduction}

Spiking Neural Networks (SNNs), named ``the third generation of neural networks'' \cite{maass1997networks}, are an emerging trend in the scientific community and they are becoming the center of attention of many research groups. A SNN operates on input spiking data (spike-encoded images, text, continuous values, etc.) and produces output spiking predictions. The spiking framework makes the implementation and interpretation of SNNs much different from their Artificial Neural Networks (ANNs) counterparts. The incentive for exploring the use of SNNs is that the human brain, arguably the most sophisticated learning machine in existence, uses much less power than popular hardware used for machine learning (Core and Graphics Processing Units, CPUs and GPUs). The human brain stands as a good example for learning algorithms that are both more accurate, as well as efficient. However, creating such biologically plausible learning algorithms is a challenging task. Many authors are investigating this topic \cite{cronin1987mathematical,guo2020self,song2000competitive,van2000stable}, but in most cases the methods developed are not fully satisfactory yet in terms of the accuracy/efficiency trade-off.

On the hardware side, neuromorphic processing units \cite{mcdonnell2014engineering,rueckauer2022nxtf} have been developed to perform computations using spike trains, similar to what happens in the brain. These chips can be analog or digital (or mixed), and consist of many small units that carry a signal, either 0 or 1. These small units mimic the neurons in the brain, and the chips perform Very Large Scale Integration (VLSI) to perform computations using these units. A layer of logical gates allows the chips to control these units, and using these one can implement algorithms on neuromorphic chips. SNNs are meant to run on neuromorphic chips in order to advance towards the full potential of biologically inspired machine learning. In this work, we focus on the algorithmic side. The method proposed herein can, in principle, be implemented to run on a neuromorphic chip with much higher efficiency.

SNNs are often used for classification, since the output of a SNN is usually a few neurons that have either 0 or 1 values, which resembles a one-hot encoded vector in classification tasks \cite{diehl2015unsupervised,lee2016training,tavanaei2019deep}. On the other hand, in Scientific Machine Learning (SciML), we are interested in using SNNs for regression, which includes approximating continuous functions, solutions of ODEs and PDEs, and operators. In \cite{spiking_deeponet}, we presented the first SNN-based method involving operator regression \cite{lu2021learning} to approximate functions and solutions of PDEs. In this paper, we aim to extend some of these ideas to time-dependent problems using a new approach.

To approximate solutions of time-dependent differential equations, one usually employs numerical methods  \cite{hairer1993solving,wanner1996solving,iserles2009first}. There is an obvious trade-off between the accuracy and the efficiency of a method, since higher accuracy requires more computational resources. Most numerical methods can be classified as either implicit or explicit. Explicit methods approximate the solution at the next timestep from the solution at the current timestep, while implicit methods obtain the approximate solution by solving an equation involving the solution at both the current and the next timesteps. %Implicit methods are usually less efficient but more accurate while explicit methods are usually more efficient but less accurate. NOT TRUE! %
All explicit methods are conditionally stable, meaning that if the stability condition is not met, the approximate solution accumulates catastrophically large errors even within a few timesteps.
In the context of time-dependent PDEs, popular explicit numerical schemes are, e.g., those based on finite differences. The idea is to define a spatio-temporal mesh, and by  approximating the derivatives of the solution on that mesh using values at the current (and previous) timestep, one can create an explicit marching scheme and use it to advance an approximation of the PDE solution over time. Explicit schemes are usually straightforward in terms of interpretation and implementation. However, they are conditionally stable, for example, for advection-dominated problems they are characterized by the Courant-Friedrichs-Lewy (CFL) condition \cite{CFL} while for the diffusion-dominated problems they are characterized by the Diffusion condition. The CFL and Diffusion conditions couple the spatial and temporal discretizations. In particular, finer spatial discretization has to be matched by a proportionally finer temporal discretization. This can become a serious bottleneck for the simulation of phenomena that require long time integration. 

%%%%%%%%%%%%%%%%%%%%%%%%%%%%%%%%%
\begin{figure}[h!]
    \centering
    \includegraphics[width=15cm]{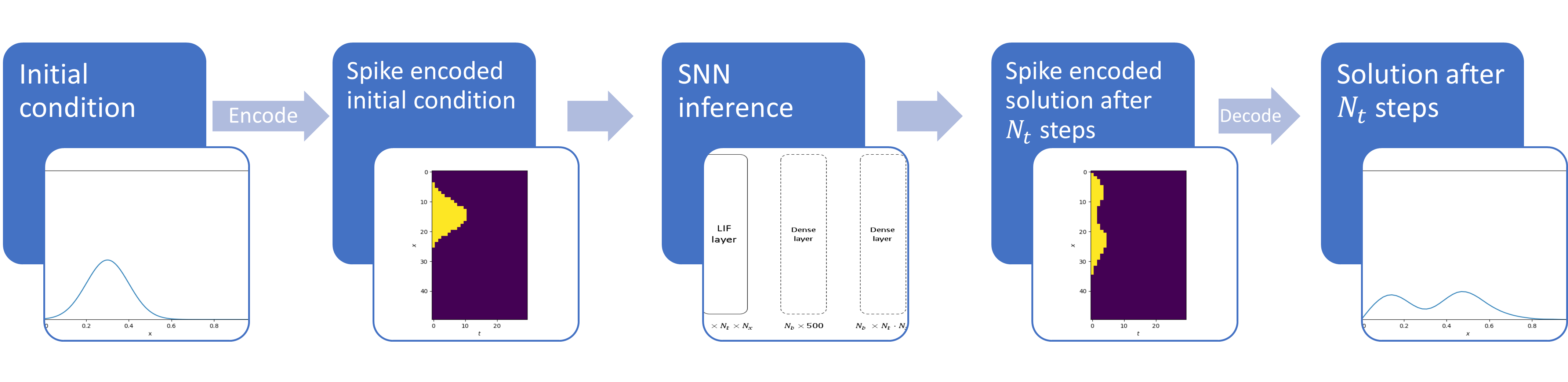}
    \caption{Illustration of the proposed approach. First, we  encode the initial condition into a spike train. Then, we feed the spiking data into a trained SNN and use it iteratively to advance the solution for many time steps. Finally, we decode the output to get the non-spiking solution.}
    \label{fig:approach}
\end{figure}
%%%%%%%%%%%%%%%%%%%%%%%%%%%%%%%%%

In this work, we develop a SNN based method to construct numerical integrators for ODEs and PDEs. In the proposed construction, the trained SNN acts as an explicit method. It receives as input an approximate solution at a previous timestep and predicts the approximate solution at the next step. The input and output of the SNN are in the form of spikes, and are able to utilize the efficiency advantages of SNNs. To this end, we encode the continuous inputs into spiking data, and after advancing them in time through the use of the SNN, we decode the output from spikes back into function values. After the network has been trained, one encodes an initial condition into spikes, applies the SNN recursively to advance the solution for many timesteps with low computational cost, thus approximating the solution for long times with high accuracy (if fine temporal discretization is used), and decodes the output. This process is illustrated in Fig. \ref{fig:approach}. Through numerical experiments, we show that the proposed approach has no stability constraints and even when not complying with the CFL or Diffusion conditions it is still able to make accurate predictions of the solution for long times. The full potential of the proposed approach in terms of efficiency requires implementation on a neuromorphic chip, which we have not pursued here.  

%% file: background.tex
\section{Related work}

In recent years, biologically plausible machine learning (ML) methods have been  attracting increasing interest. The three popular directions are: a) converting a pre-trained artificial neural network (ANN) to a SNN \cite{diehl2016conversion,deng2021optimal}; b) building surrogate models so that back propagation can be implemented in SNNs \cite{neftci2019surrogate,xiang2022spiking}; and c) using a pure SNN with only biologically plausible components \cite{diehl2015unsupervised}. In \cite{javanshir2022advancements}, the authors present a review of a large variety of SNN based methods and compare them. Since most SNN architectures are designed for classification, specifically classifying handwritten digits \cite{diehl2015unsupervised}, they can be compared in terms of accuracy and efficiency on this specific task. There has also been work on the problem of classifying images \cite{zheng2021going,rathi2020enabling}.

With respect to the three popular directions mentioned above, each has advantages and disadvantages. The first is about converting a trained ANN into an SNN. Two drawbacks of this approach are that it is hard to convert complex ANN architectures into SNNs (for example, supporting all kinds of activation functions such as hyperbolic tangent, sinusoidal, etc.), and also training the ANN still requires a lot of  computational effort. The second approach discusses building surrogate models to train a SNN. These methods are as accurate and efficient as the developed surrogates. Lastly, there are SNNs that use Hebbian learning rules \cite{song2000competitive} for training, which however are highly unstable. Even if one manages to stabilize them and train them correctly, they currently underperform in terms of accuracy, and cannot compete with current state-of-the-art ANNs.

Scientific Machine Learning (SciML) is becoming increasingly popular for problems in computational science and engineering. For example,  Physics-Informed Neural Networks (PINNs) \cite{raissi2019physics} have been applied to various physical and biological systems with considerable success. However, like other learning methods, they consume a lot of computational resources and require a lot of data for training, especially for complex problems (high dimensions, high frequencies, etc.). There are some works on computing forward gradients \cite{forward_prop} and eliminating the need for back-propagation \cite{rumelhart1986learning}, which is the most expensive part of training the neural network. 
%Another drawback of PINNs is extrapolation. PINNs are usually trained in a specific domain, and studies show their ability to correctly interpolate data to all points (continuously) in that domain. However, extrapolating to points outside that domain is still a challenging task for PINNs and many other ODE and PDE solvers.

We began exploring this field with the goal of using SNNs to enhance the training of machine learning based methods for solving differential equations. The first step was to use SNNs for regression \cite{spiking_deeponet}. We proposed a method to encode input numeric data ($x$) and train a SNN to infer the values of a function ($y(x)$). We also proposed using a deep neural operator, DeepONet \cite{lu2021learning}, for the same function regression task, essentially approximating a spiking (and not continuous) operator. We used back-propagation for training the networks, as a surrogate for a synapse. The proposed method operates on structured grids, and while we managed to show interpolation capabilities, extrapolation was not investigated in that work.

Time marching schemes are essential in solving time-dependent PDEs. They have been developed for a long time to achieve improved efficiency and accuracy \cite{iserles2009first}. Explicit schemes are usually faster, but less accurate and only conditionally stable, meaning that without satisfying the CFL or Diffusion conditions the solution may diverge after just a few time integration steps. On the other hand, implicit schemes  are often more accurate and without timestep constraints (up to second-order), but require a much larger computational effort. There are also Implicit-Explicit schemes (IMEX) that combine both approaches to achieve better performance \cite{ascher1995implicit}.

ML based methods for long time integration have also been proposed. Properly trained, PINNs have shown an ability to extrapolate for longer times \cite{meng2020ppinn,wang2021long}. In addition, DeepONets can be trained for this task as well \cite{bubbles}. In \cite{ovadia2021beyond} the authors proposed a ML based explicit scheme that is stable even when violating the stability conditions. A noticeable drawback of such methods is that they take long time to train, and the PINN based methods require training for each problem being investigated. A less noticeable drawback is that machine learning methods often reach a limit of accuracy. Whereas numerical solvers may achieve an error norm of, e.g., $10^{-10}$ per step, ML based methods may only achieve at best a $10^{-5}$ error norm, resulting in large error accumulation over time. As a result, they are able to forecast accurately only for relatively short time intervals.

On the hardware side, neuromorphic chips have been developed for fast implementation of spiking neural networks. (e.g., Intel's Loihi 2 \cite{davies2018loihi, orchard2021efficient}). These chips are using either an analog or a digital interface, and their architecture consists of many small units that can either have a current flowing through or be idle (spiking or not). Software packages have been developed to interact with such chips, defining the logical gates that switch the currents between the tiny units, enabling the users to train SNNs. In addition, simulators (such as LAVA \cite{lava} for Loihi 2) have been made available, enabling researchers to develop code that can be implemented on the neuromorphic chip, while testing it first on non-neuromorphic machines.

%% file: sms.tex
\section{Spiking Marching Scheme}

\subsection{Formulation as a data-driven problem}\label{section:data_driven}

We propose using a machine learning based method, which utilizes SNNs for the time integration of ODEs and PDEs. We employ supervised learning to train the neural networks. To do that, we create a dataset of samples consisting of inputs and outputs. The goal is to create a SNN based numerical scheme, where the input is the solution at the current (and if needed previous) timestep(s), and the output is the solution at the next time step. In the case of ODEs, the solution at the current timestep is specified by the state vector of the system. In the case of PDEs, the solution at the current timestep is specified on a pre-defined spatial grid. 

Without loss of generality we can define a PDE as follows:

\begin{equation}
    \begin{aligned}
        \mathcal{L}_{\textbf{x},t}(u;k)&=f, \quad \textbf{x} \in \Omega, \quad t\in[0, T] \\
        \mathcal{IC}_{\textbf{x},0}(u)&=u_0, \quad \textbf{x} \in \Omega \\
        \mathcal{B}_{\textbf{x},t}(u)&=u_b, \quad \textbf{x} \in \partial\Omega, \quad t\in[0, T]
    \end{aligned} \qquad,
    \label{pde}
\end{equation}
where $\textbf{x}$ is the spatial coordinate, $t$ is the temporal coordinate, $\Omega$ is a domain in $\mathbb{R}^d$. The terms $\mathcal{L}_{\textbf{x},t}$, $\mathcal{IC}_{\textbf{x},0}$ and $\mathcal{B}_{\textbf{x},t}$ represent the partial differential operator, the initial condition, and the boundary operator, respectively. The quantities $k(\textbf{x})$ reflect parameters that the differential operator can depend on, and $f$, $u_0$ and $u_b$ are the forcing term, the initial condition and boundary condition, respectively. The quantity $u(\textbf{x},t)$ is the solution of the PDE. We present a general formulation of differential operators since the method proposed in this paper can be applied to a variety of problems of this formulation with few adjustments only. In this work, we focus on a small set of problems with specific parameters, initial and boundary conditions. The exact problem parameters are given in their specific subsections of Sec. \ref{results}. We note that an ODE can be defined as a special case of \eqref{pde} by focusing only on one independent variable (either the spatial or temporal).

We discretize the problem in space and time (for simplicity, we use a uniform grid). In the case of one spatial dimension and one temporal dimension we use the following notation: $u^n_i=u(x_0+i\Delta x; n\Delta t)$, $i=0, ..., N_x$ and $n=0, ..., N_t$, where $N_x$ and $N_t$ are the number of points in the spatial and temporal discretizations, respectively. Given this notation, a couple of examples of finite difference approximations (of second-order accuracy in $\Delta t$) for the first and second temporal derivatives are:
\begin{equation}\label{fd_derivatives}
    \begin{aligned}
        u_{t} &\approx \frac{u^{n+1}_i-u^{n-1}_i}{2\Delta t} \\
        u_{tt} &\approx \frac{u^{n+1}_i-2u^{n}_i+u^{n-1}_i}{\Delta t^2}
    \end{aligned}\qquad ,
\end{equation}
with similar expressions for spatial derivatives. 
% Higher order derivatives are also possible to compute in this fashion using higher order Runge-Kutta methods \cite{} and various other approaches, but for the purpose of this work, the first and second order derivatives are sufficient.
For example, if the PDE contains up to second-order temporal derivatives, we can create an explicit time marching scheme by using (\ref{fd_derivatives}) to re-write the approximation of the PDE as $u^{n+1} = \mathcal{L}_*(u^{n}, u^{n-1})$ (for ease of notation we have omitted the spatial index). The expression $\mathcal{L}_*$ is often referred to as the scheme or a kernel, and its exact form depends on the order of the method and the form of the PDE. Using the initial conditions, given usually as $u^0, u^1$, one can compute $u^2$ and iteratively continue until the terminal timestep number $N_t$ is reached.

To formulate this process as a data-driven problem, a possible approach is to use $u^{n}, u^{n-1}$ as input data, and train a learning algorithm to infer $u^{n+1}$ as the output data. In this work, we also consider using only $u^{n}$ to infer $u^{n+1}$. In addition, we also tried marching with larger timesteps by using $u^{mn}$ to infer $u^{m(n+1)}$, where $m$ is a positive integer number larger than one. The data used to train the SNN based model were created using appropriately accurate existing time marching schemes for ODEs and PDEs. 

\subsection{SNN Formulation}

\subsubsection{Spike encoding-decoding}\label{section:encoders}

\begin{figure}
     \centering
     \begin{subfigure}[h]{0.49\textwidth}
         \centering
         \includegraphics[width=\textwidth]{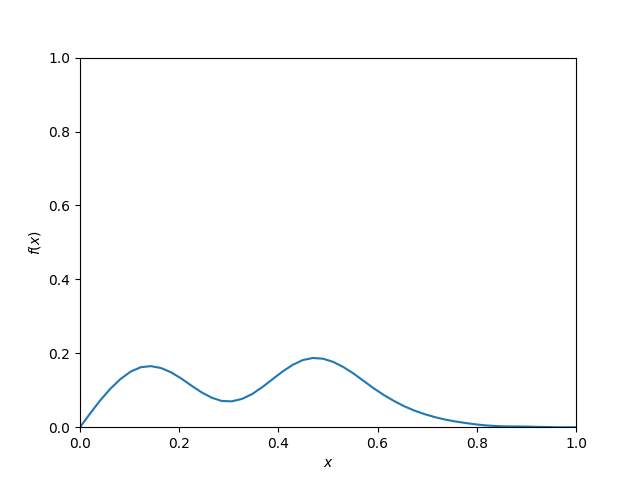}
         \caption{A sample function $f(x)$.}
     \end{subfigure}
     \hfill
     \begin{subfigure}[h]{0.5\textwidth}
         \centering
         \includegraphics[width=\textwidth]{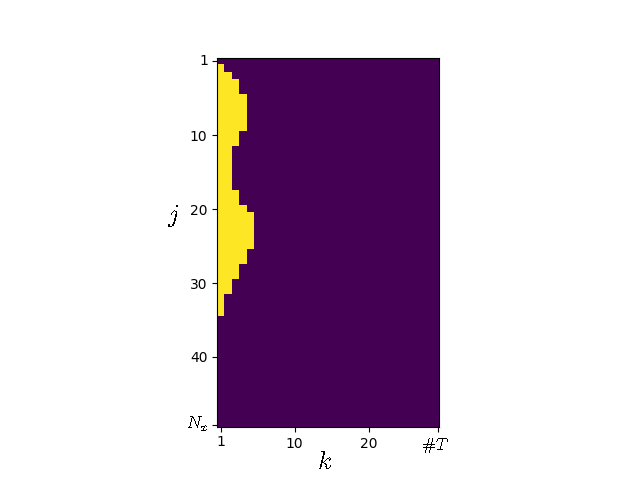}
         \caption{The sample function $f(x)$ encoded into spikes. Pixels with value 1 are in yellow and pixels with value 0 are in dark blue. The vertical axis corresponds to the spatial discretization ($j=1, ..., N_x$) and the horizontal axis corresponds to the timesteps in the spikes train ($k=1, ...,\#T$).}
     \end{subfigure}
     \hfill
     \caption{Spike encoding illustration.}
     \label{fig:overall_encoding}
\end{figure}

We encode the input data into spikes. To encode the data we introduce a new axis for the temporal behavior of the spikes, and the continuous function values are encoded into binary spikes. We use an approach inspired by the lower triangular encoding presented in \cite{spiking_deeponet}. Suppose that we want to encode a function $f(x)$ which is defined on the spatial interval $[0,1].$ To encode the function, first, we consider a discretization of the spatial axis. Second, we define a temporal axis and we use the number of temporal nodes chosen for the new temporal axis to discretize the function values. The upper and lower limits of this discretization are chosen based on the function values and may very between functions. Then, we add spikes based on the function value. For example, if the number of nodes in the temporal axis is $10$ and the function value at a point $x$ is $f(x)=0.4$ (assuming the function values are between $0$ and $1$), the spike encoding for that value is $(1111000000)$, so it has four ones as the spikes. We do this for every value of the function and get a two dimensional array, where the first dimension is the location on the spatial axis and the second is the location on the temporal axis. The data is binary as desired. In this way we exploit the sparsity benefit of SNNs. The drawback of this encoding is that we lose some information. After encoding, we cannot reconstruct the original signal, and the loss of information depends on the discretization of the temporal axis. A more detailed description of the  encoding mechanism can be found in \cite{spiking_deeponet}. An example of encoding a function is given in Fig. \ref{fig:overall_encoding}. The function is encoded using 30 steps in the spikes train (the horizontal axis of Fig. \ref{fig:overall_encoding}(b)). A similar approach is applied for the cases presented in Sec. \ref{results}.

In some of the numerical experiments we also consider other encoding methods. There are two reasons for this: a) we would like to show that the spiking marching scheme can operate with other encoding methods, in other words, it is not limited to the aforementioned encoding method, and b) we are interested in the effect of different encoding methods on the performance of the spiking marching scheme. The second encoding method we use is the popular rate encoding approach, where one normalizes the data points (values of the function) into the $[0, 1]$ interval and for each value, creates Poisson trials with expected value set as the normalized function value. The number of trials is set to the number of time steps of the spikes train, which is denoted by $\#T$. In that way, the function representation is expanded into a $N_x\times \#T$ binary images. The third encoding method we use is representing the number using the float32 representation of it, which is indeed binary and in this case the number of steps in the spikes train is set to 32. For more information about these encoding methods we refer the readers to \cite{spiking_deeponet}. To decode the data (at the end of the process) we use the exact opposite of the encoding procedures mentioned here. 

We compare the performance of the spiking marching scheme with lower-triangular encoding to that with either rate encoding or the float32 encoding. We point out which encoding method was used for each of the experimental setups in Sec. \ref{results}. The main motivation is showing that the spiking marching scheme is not encoding-specific, and can work with different types of encoding methods. Also, we are interested in choosing the encoding method that gives the best overall performance of the spiking marching scheme (most accurate and efficient).

Encoding methods induce a certain error for the computation. When encoding the signal, we use a fixed number of spike train steps ($\#T$), and the smaller this number is, the more information is lost. If one wishes to reconstruct the encoded signal accurately, a large $\#T$ should be used which then requires more computational resources to train the network. Facing this trade-off between efficiency and accuracy, a fixed $\#T$ is chosen for each problem in the numerical tests. We investigate the impact of the encoding error by taking a reference signal, encoding it using each of the three encoding methods, decoding it back and computing the $L_2$ error. The reference signal taken is $\eta_{ref}(x_j) = e^{-\frac{1}{2}}\left(\frac{x_j-\mu}{\sigma}\right)^2, x_j\in\{j\Delta x\}_{j=0}^{N_x}$ (linearly spaced $[0,1]$ interval with $\Delta x = 0.01$), with $\mu=\frac{1}{2}$ and $\sigma=0.1$. The error is computed by: $||\eta_{encoded}(x)-\eta_{ref}(x)||_{L_2}=\frac{1}{N_x}\sum_{j=1}^{N_x}|\eta_{encoded}(x_j)-\eta_{ref}(x_j)|^2$. The error comparison is given in Table \ref{tab:encoding_errors}. We can observe that the error for the lower triangular encoding is much lower than the error for the rate encoding methods. The float32 encoding error is, as expected, the lowest and in fact it is at the level of machine zero. However, it does not necessarily translate to better performance of the spiking marching scheme. We explore in Sec. \ref{results} the impact of the encoding methods on the prediction errors for various systems.

\begin{table}[h]
    \centering
    \begin{tabular}{|c|cc|}
        \hline
        & $L_2$ & Relative $L_2$\\
        \hline
        Lower Triangular & 2.1195e-05 & 3.3822e-04 \\
        Rate & 5.8053e-03 & 9.2640e-02 \\
        Float32 & 2.8336e-17 & 4.5218e-16 \\\hline
    \end{tabular}
    \caption{Encoding errors comparison. The $L_2$ and relative $L_2$ errors of encoding and decoding a reference signal using the three encoding methods employed in this study (see text for details).}
    \label{tab:encoding_errors}
\end{table}

The quality of the decoding (sharpness of the recovered signal) is also determined by the number of time steps we choose to encode the data. There is a trade-off between the accuracy (finer spike discretization is more accurate) and efficiency (more spikes require more computations). We experimented with several values of spike numbers, eventually choosing one that is suitable for each problem by trial-and-error.

\subsubsection{Network architecture}

\begin{figure}[h!]
    \centering
    \includegraphics[width=15cm]{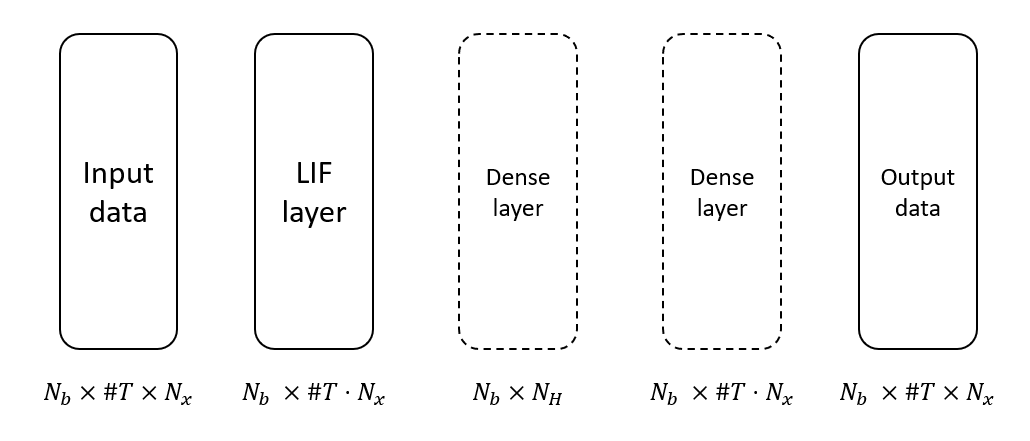}
    \caption{Schematic of the network architecture, consisting of an input layer followed by a LIF layer. Then, we flatten the output and use two dense layers that are trainable (marked in dashed lines). The output layer has the same size as the input layer.}
    \label{fig:arch}
\end{figure}

For the SNN architecture we use a relatively simple model, with a Leaky Integrate and Fire (LIF) layer followed by two dense layers. We use $\#T$ steps for number of time steps of the SNN. Note that $\#T$ is different from the number of timesteps of the physical problem $N_t$. The former regards the number of steps in the spikes train (see a detailed explanation in \cite{spiking_deeponet}). $N_x$ is the number of points in the spatial discretization, and the first dense layer is of size $N_H$ (chosen equal to $500$ for all of the numerical experiments in the current work). The second dense layer has the same size ($\#T \cdot N_x$) as the output layer. Both dense layers are followed by a Rectified Linear Unit (ReLU) activation. Note that we flatten the inputs (noted by $\#T\times N_x$ and $\#T\cdot N_x$) in order to apply the intermediate dense layers. $N_b$ is the size of the batch (number of samples). A schematic of the architecture used in this work is given in Fig. \ref{fig:arch}. The loss function we use for training the model is the Sigmoid Cross Entropy (SCE) loss \cite{cui2019class}.

The approach proposed  here is proof-of-concept for the use of SNNs to construct numerical integrators for differential equations. Thus, we are fully aware that the simple architecture we propose can be improved. Some studies use an architecture that mixes LIF layers with dense layers. For example, one can start with a dense layer, followed by a LIF layer, followed by a dense layer, followed by another LIF layer (see e.g. tutorials of SNNTorch, a popular package for training SNNs \cite{eshraghian2021training}). A different architecture is proposed in the tutorial of the LAVA \cite{lava} framework, in particular an architecture consisting of  three successive (LIF layer - dense layer) pairs. Because the LIF layer has a thresholding operation inside of it, performing back-propagation  requires computing gradients through the use of a surrogate model, which can introduce errors. We chose to use a LIF layer followed by two dense layers because it does not require a surrogate model for the computation of back-propagation gradients (the thresholding is done before the dense layers), and is sufficiently expressive for presenting the capabilities of the proposed method.

\subsection{Accuracy metrics}

To assess the accuracy of the proposed approach, we use two metrics. The first is the error committed during a one timestep prediction (inference). We compute the $L_2$ error of the estimate of the solution at a certain timestep given by the proposed approach with respect to the reference solution, given as input the reference solution values. Specifically, we use the reference solution $u^{mn}$ (or both $u^{m(n-1)},u^{mn}$) as input to the SNN to estimate $u^{m(n+1)}$. This is what the network is trained to do and the errors are expected to be relatively small. This evaluation serves as a sanity check, and if it fails it signifies a problem with the training of the model.

The second metric is what we call cascade prediction (inference). We start from the solution at the initial timestep and predict the solution at the next timestep. We then use the predicted solution timestep ({\it not} the reference solution like before) to predict the solution at the next timestep, and so on. Since every prediction  commits an error, the  errors can accumulate and cause divergence of the predicted solution from the reference one. This is a more challenging accuracy metric, and is the main criterion used to monitor the accuracy of the proposed approach. 

For both the one-step and cascade metrics, we are interested in observing the ability of the spiking marching scheme to make correct predictions: a) for times {\it inside} the temporal interval used for training (interpolation), and b) for times {\it outside} the temporal interval the network was trained on (extrapolation). All marching schemes suffer from accumulating numerical errors over time, and seek to keep the cascade inference error from growing too fast. 

%% file: results.tex
\section{Results}\label{results}

\subsection{ODE Examples}

We have applied the proposed approach to solve the initial-value problem for two ODE systems. The first is the Van der Pol oscillator \cite{guckenheimer1980dynamics} while the second is the Lorenz system, which is chaotic and thus presents a challenge for any time marching scheme \cite{lorenz1963deterministic}.

\subsubsection{Van der Pol oscillator}

The Van der Pol equation is a mathematical model for the motion of an oscillator, which is nonlinearly damped. It is a second-order ODE, which can be rewritten as a system of two first order ODEs:
\begin{align*}    \frac{dx}{dt} &= y \\
        \frac{dy}{dt} &= \mu (1 - x^2)y - x
\end{align*}
where $x$ and $y$ are the position and velocity of the oscillator, respectively, $t$ is the temporal variable, and $\mu$ is a parameter. We solve the Van der Pol system in the interval $[a,b]$. For this test, we choose $\mu = 2, a = 0, b = 10$ with initial condition $x(a)=1$ and $y(a)=0.$ 

\begin{figure}
     \centering
     \begin{subfigure}[b]{0.49\textwidth}
         \centering
         \includegraphics[width=\textwidth]{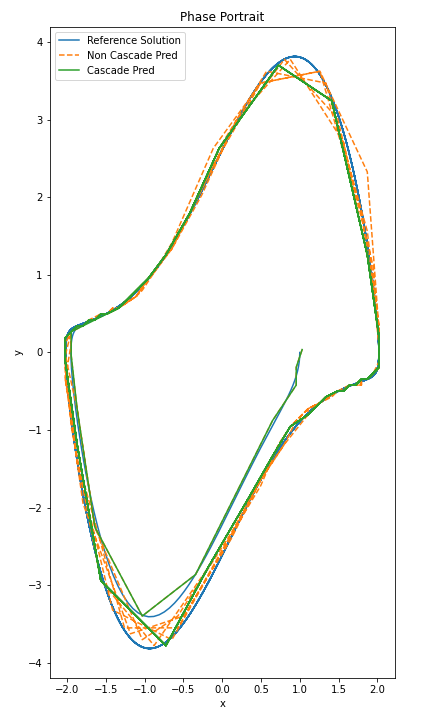}
         \caption{Trajectory prediction.}
         \label{fig:vanderpol_trajectory}
     \end{subfigure}
     \hfill
     \begin{subfigure}[b]{0.5\textwidth}
         \centering
         \includegraphics[width=\textwidth]{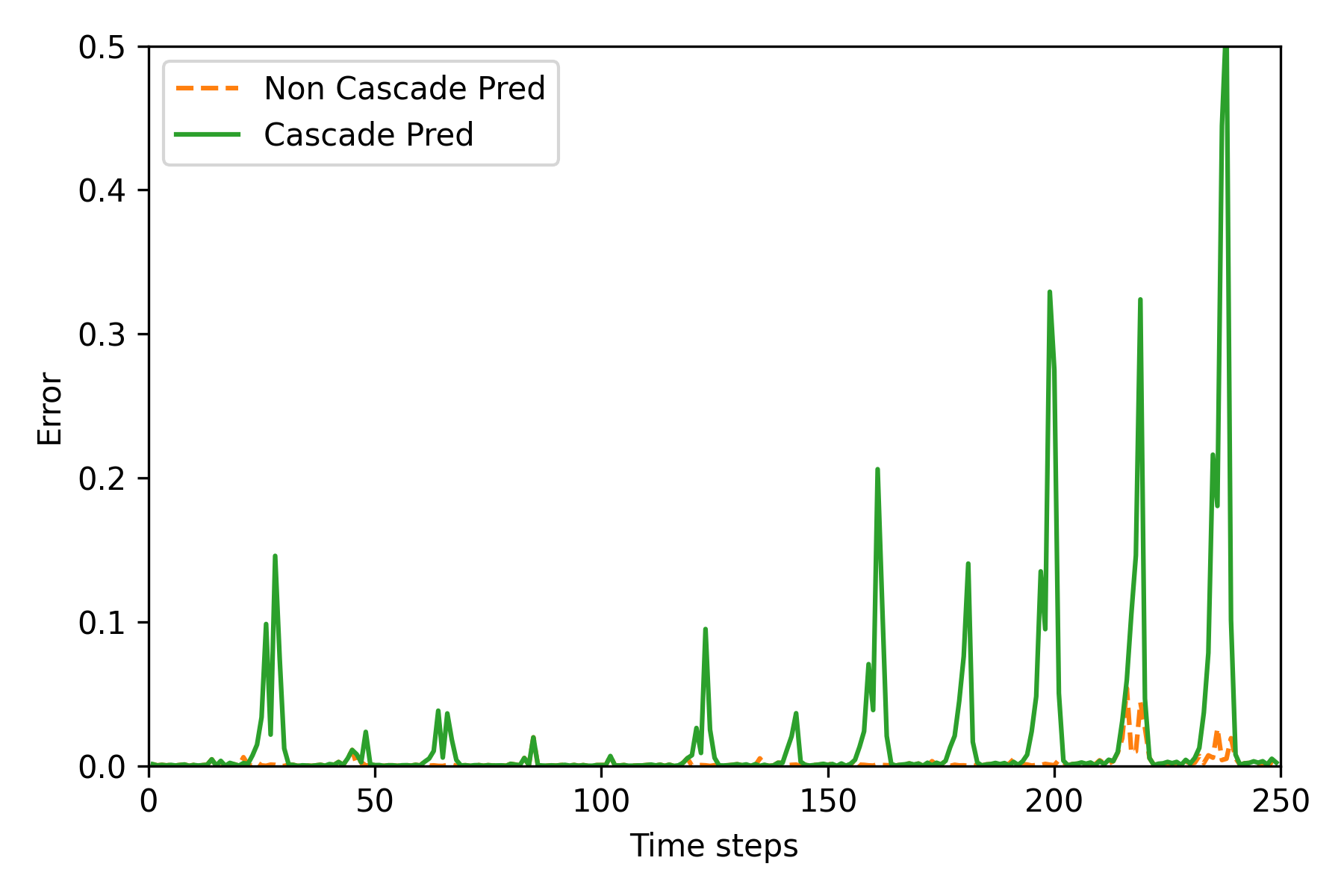}
         \caption{Error evolution over timesteps. The first 200 timesteps show interpolation and the remaining 50 show extrapolation. Here the error at timestep $n$ is defined as $\frac{1}{2}((x^n-\Tilde{x}^n)^2+(y^n-\Tilde{y}^n)^2)$, where $x, y$ are the prediction and $\Tilde{x}, \Tilde{y}$ are the reference solution.}
         \label{fig:vanderpol_error}
     \end{subfigure}
     \hfill
     \caption{Van der Pol oscillator. Encoding method: Lower triangular encoding.}
     \label{fig:vanderpol}
\end{figure}

To create the dataset for training we use a numerical initial value problem ODE integrator in Python (from the SciPy library \cite{2020SciPy-NMeth}). We compute the solution for $N_t = 2,500$ steps and sub-sample every 10th step. This means that in the new dataset, we effectively use a $10$ times larger $\Delta t$, ending up with $250$ steps for training and testing. The first 200 steps are used for training (interpolation) and the last 50 for testing (extrapolation). We use both previous and current time steps to infer the next time step, so each sample is of size $4$ (the $x, y$ components of both previous and current steps), which substitutes in the ODE setting the notion of spatial discretization  (see $N_x$ in Fig. \ref{fig:arch}). We then encode the data using the lower triangular encoding with $\#T = 100$ spike train steps. The model is trained with $5000$ epochs. Fig. \ref{fig:vanderpol_trajectory} shows the trajectory prediction using SMS for $N_t = 250$ timesteps with both cascade and non-cascade predictions. The evolution of the error with time is given in Fig. \ref{fig:vanderpol_error}.

We observe the highly accurate one-step and cascade predictions in the first 100 steps (within the interpolation interval). Also, the error accumulates slowly outside of the training interval (beyond 200 steps) where the marching scheme performs extrapolation. We note that the accuracy of the trained model is promising, given that it takes steps which are 10 times larger than the ODE integrator used to produce the reference solution.

\subsubsection{Lorenz system}\label{section:lorenz}

\begin{figure}[!ht]
    \centering
    \begin{subfigure}[b]{0.49\textwidth}
         \centering
         \includegraphics[width=0.8\textwidth]{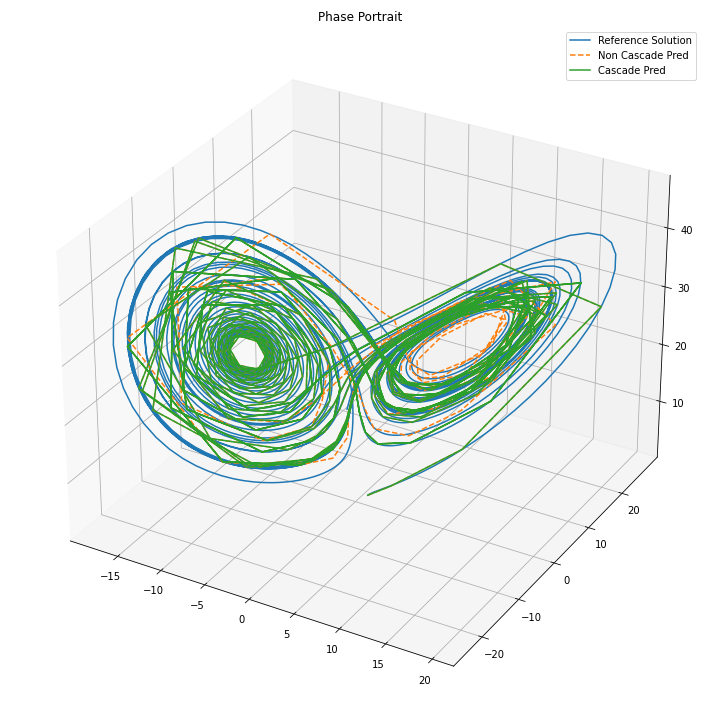}
         \caption{Trajectory prediction, lower triangular encoding.}
     \end{subfigure}
     \hfill
     \begin{subfigure}[b]{0.5\textwidth}
         \centering
         \includegraphics[width=0.8\textwidth]{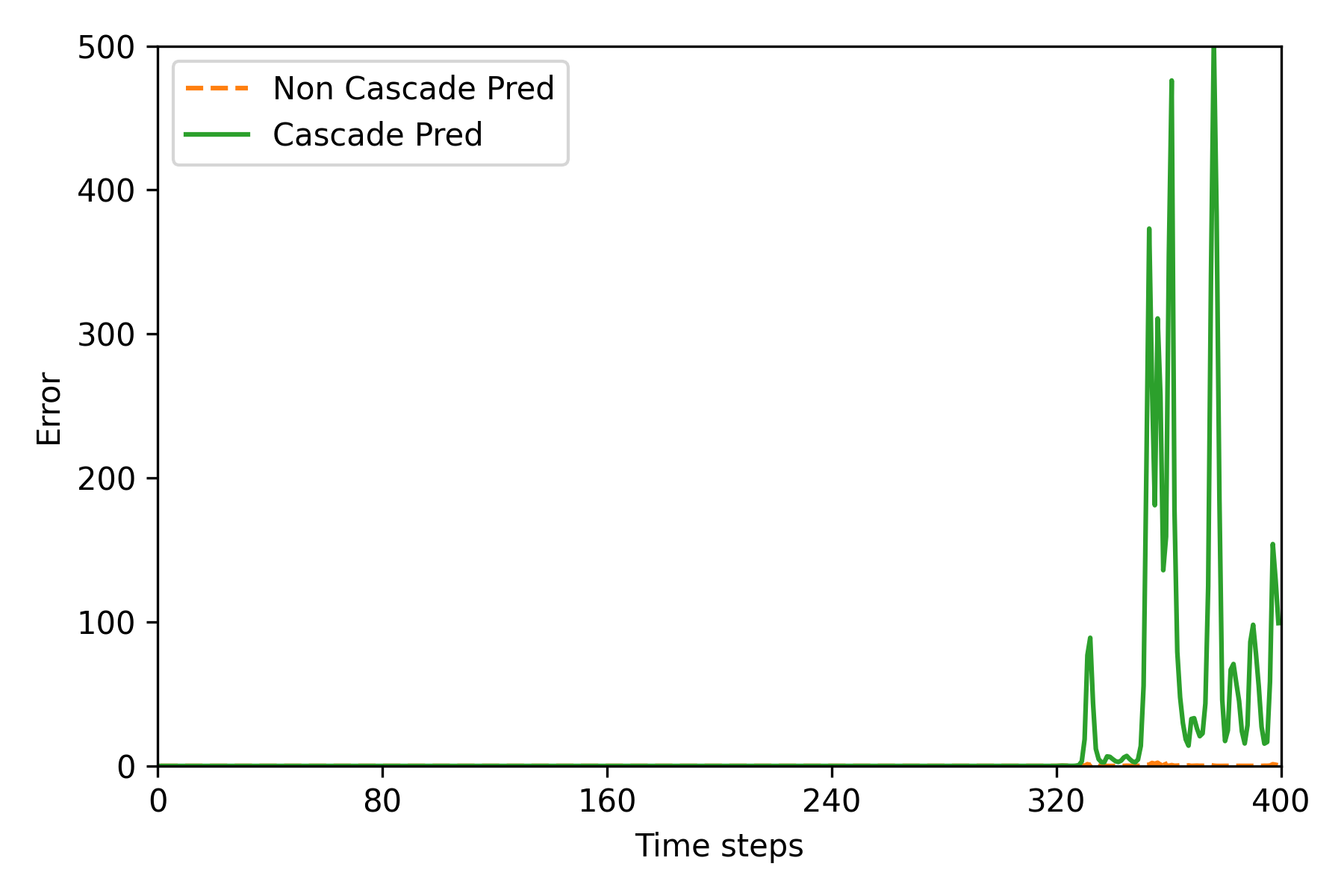}
         \caption{Evolution of error over timesteps. The first 320 timesteps show interpolation and the remaining 80 show extrapolation. Here the error at timestep $n$ is defined as $\frac{1}{3}((x^n-\Tilde{x}^n)^2+(y^n-\Tilde{y}^n)^2+(z^n-\Tilde{z}^n)^2)$, where $x, y, z$ are the prediction and $\Tilde{x}, \Tilde{y}, \Tilde{z}$ are the reference solution.}
     \end{subfigure}
     \hfill
     \begin{subfigure}[b]{0.49\textwidth}
         \centering
         \includegraphics[width=0.8\textwidth]{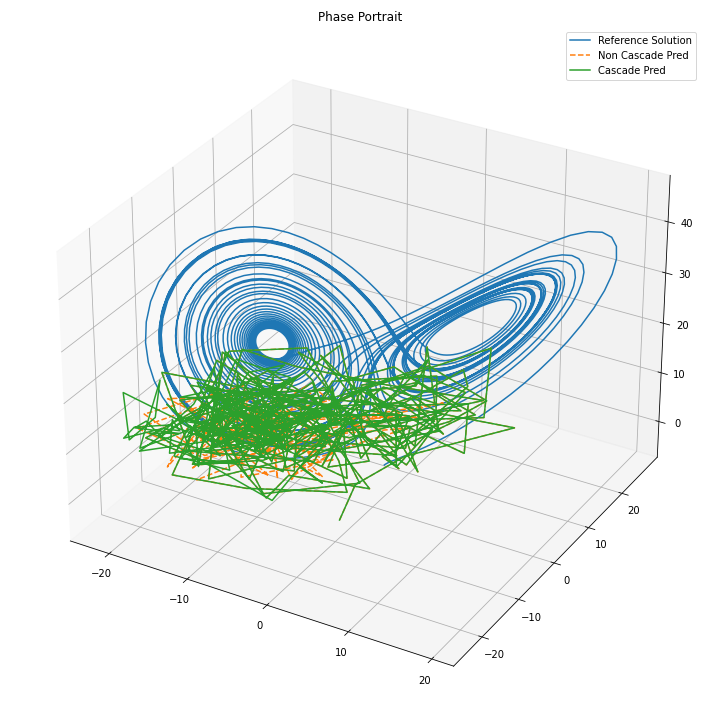}
         \caption{Trajectory prediction, rate encoding.}
     \end{subfigure}
     \hfill
     \begin{subfigure}[b]{0.5\textwidth}
         \centering
         \includegraphics[width=0.8\textwidth]{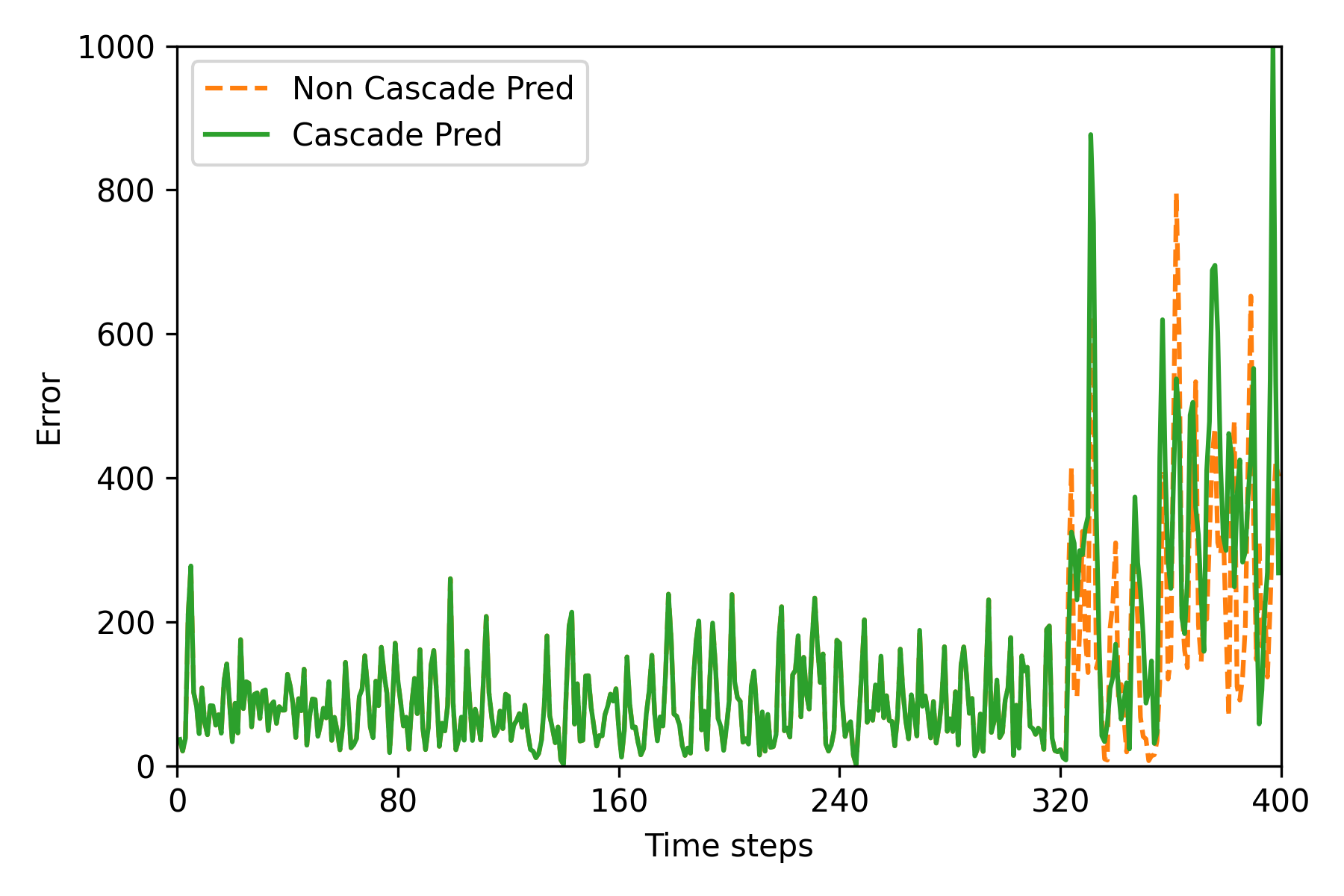}
         \caption{Similar configuration as in (b) but with rate encoding.}
     \end{subfigure}
     \hfill
     \\
     \begin{subfigure}[b]{0.49\textwidth}
         \centering
         \includegraphics[width=0.8\textwidth]{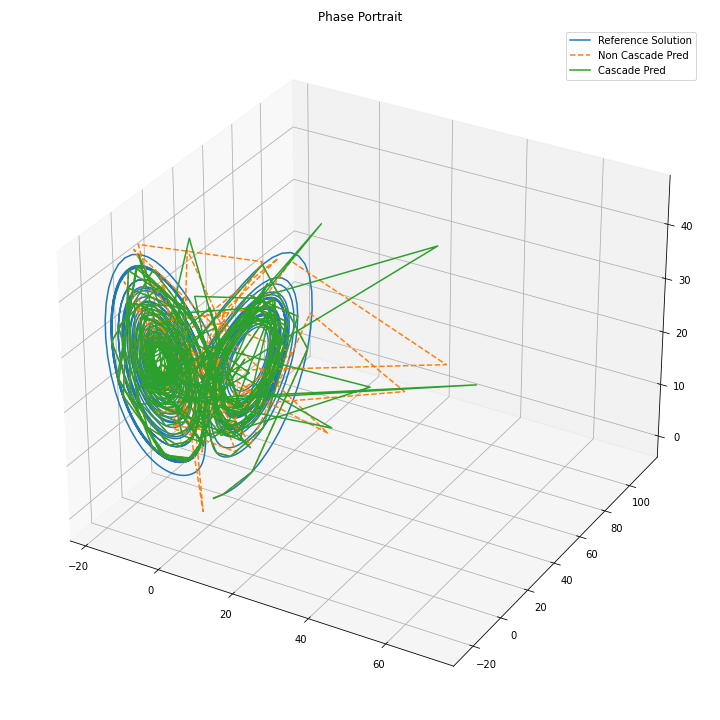}
         \caption{Trajectory prediction, float32 encoding.}
     \end{subfigure}
     \hfill
     \begin{subfigure}[b]{0.5\textwidth}
         \centering
         \includegraphics[width=0.8\textwidth]{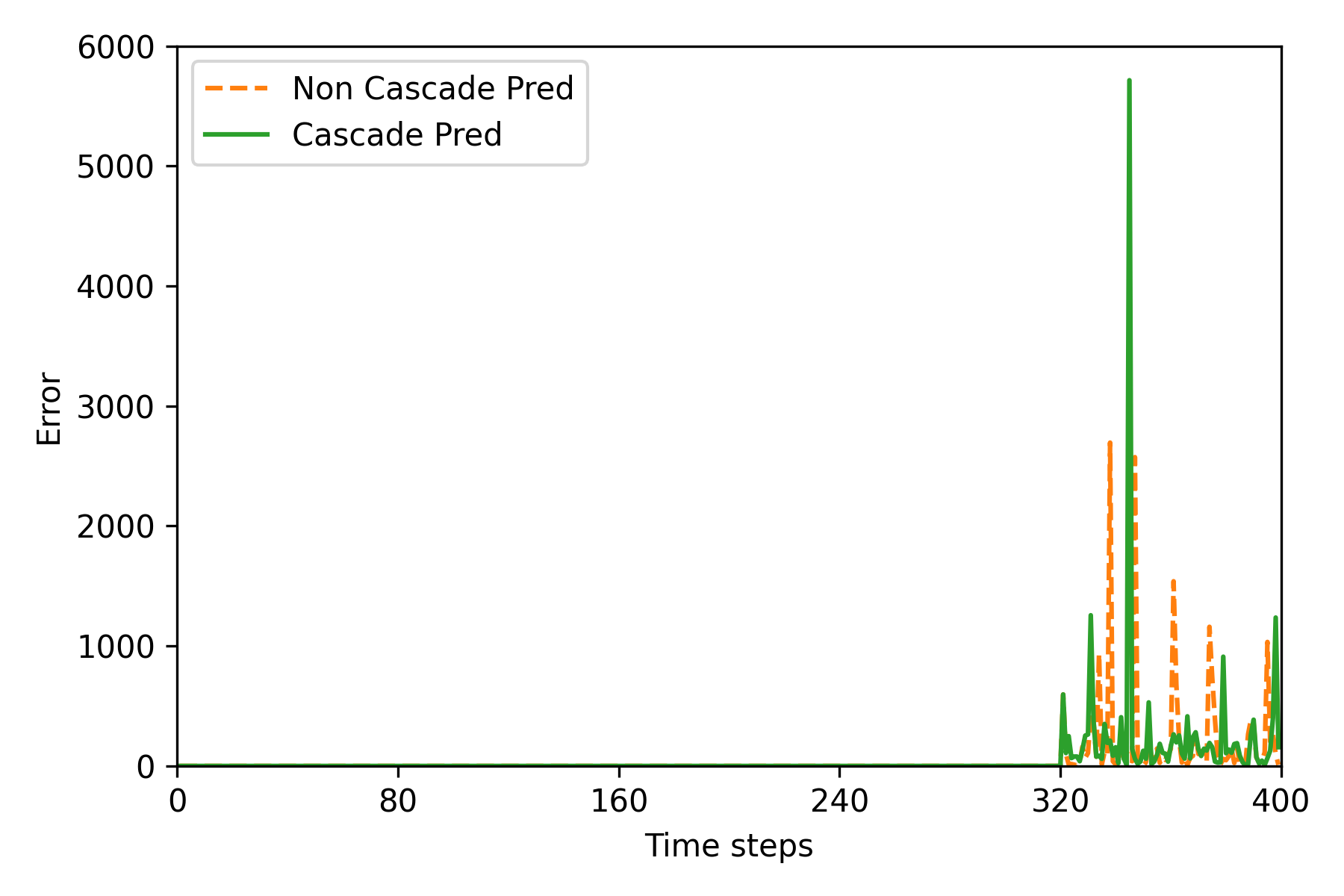}
         \caption{Similar configuration as in (b) but with float32 encoding.}
     \end{subfigure}
     \hfill
     \caption{Lorenz system. Results for various encoding schemes.}\label{fig:lorenz}
\end{figure}

The Lorenz system is given by:
\begin{align*}
 \frac{dx}{dt} &= \sigma (y - x) \\
        \frac{dy}{dt} &= x (\rho - z) - y \\
        \frac{dz}{dt} &= xy - \beta z
\end{align*}
where $x,y,z$ describe the state of the system, $t$ is the temporal variable and $\sigma, \rho$, and $\beta$ are parameters. For this experiment we choose $\sigma = 10, \rho = 28$, and $\beta = \frac{8}{3}$, which are common choices for this system and lead to chaotic behavior. The initial state is set to $(x_a, y_a, z_a)=(1, 1, 1)$ and the time interval is set to $[a,b]=[0, 40]$ with $N_t = 4,000$ timesteps. We then sub-sample the dataset and keep only 400 timesteps by using every 10th timestep in the original dataset. In this case we use both previous and current time steps to infer the next time step, so each sample is of size $6$ (the $x, y, z$ components of both previous and current steps), which substitutes in the ODE setting the notion of spatial discretization  (see $N_x$ in Fig. \ref{fig:arch}). For the Lorenz system we also experiment with using only the current time step to infer the next time step as mentioned in Sec. \ref{section:data_driven}. We achieve good results for the interpolation part, but for extrapolation the network struggles to compete with the two steps method. We encode the dataset using the lower triangular encoding, this time, with $\#T=1,000$ steps. The network is trained for $5000$ epochs on the first $320$ timesteps. To show the extrapolation capabilities, we also predict the solution for the remaining $80$ timesteps. The results are given in Fig. \ref{fig:lorenz}(a-b).

In this case, while the interpolation is accurate, extrapolation is rather difficult for the SMS. However, the error is still relatively small, and is not exploding by many orders of magnitude every iteration as often experienced with certain explicit numerical methods (demonstrated in \cite{ovadia2021beyond}).

For the Lorenz system we also demonstrate (see Fig. \ref{fig:lorenz}(c-f)) the performance of the rate and float32 encoders which were defined in \ref{section:encoders}.  We observe that for the rate encoding, the results suffer from large errors (yet not exponentially increasing in time). For the float32 encoding, we observe that the network was able to accurately infer the solution for both cascade and non-cascade predictions, but extrapolation induces an error that is about an order of magnitude larger than the one for the lower triangular encoding.

\subsection{PDE Examples}

We have also applied the proposed approach to the wave and heat equations in one spatial dimension. While these PDEs are straightforward to analyze, both have interesting features that provide insight into the performance of the SMS.

\subsubsection{Wave equation}\label{section:wave}

The one dimensional wave equation we solve is given by:
\begin{align*}
    \frac{\partial^2u}{\partial t^2} = c^2 \frac{\partial^2u}{\partial x^2}&, \qquad t\in [a, b], \quad x\in[x_{min}, x_{max}], \\
    u(0, x) = u_0(x)&, \qquad x\in[x_{min}, x_{max}], \\ 
    \frac{\partial u(0, x)}{\partial t} = v_0(x)&, \qquad x\in[x_{min}, x_{max}], \\
    u(t, x_{min}) = u(t, x_{max}) = 0&, \qquad t\in [a, b] , 
\end{align*}
where $u(t, x)$ is the wave pressure, $x$ is the spatial variable, $t$ is the temporal variable, $u_0(x)$ is the initial wave pressure, $v_0(x)$ is the initial wave velocity, chosen as $v_0(x)\equiv0$ in this experiment. The zero Dirichlet boundary condition set in the last equation causes reflective scattering from the edges of the domain (set by $x_{min}$ and $x_{max}$). For this simulation we choose $x_{min} = 0$, $x_{max} = 1$, and the spatial increments are $\Delta x = 0.01$. The wave velocity $c$ is chosen to be 1. The initial time is set to $a = 0$ and the temporal increments are chosen to satisfy the Courant-Friedrichs-Lewy condition: $\Delta t = \frac{\Delta x}{c} = 0.01$. We integrate for $N_t = 1,000$ steps so $b = \Delta t N_t = 10$, and use both previous and current time steps to predict the next time step. We do not sub-sample the dataset in this case, and choose $\#T = 100$ for encoding the data. We choose as initial condition the function:
\begin{align*}
    f(x) = e^{-\left(\frac{x - \overline{x}}{0.05}\right) ^ 2}, \quad \overline{x} = \frac{x_{max} - x_{min}}{2}
\end{align*}

\begin{figure}[!ht]
     \centering
     \begin{subfigure}[t]{0.48\textwidth}
         \centering
         \includegraphics[width=\textwidth]{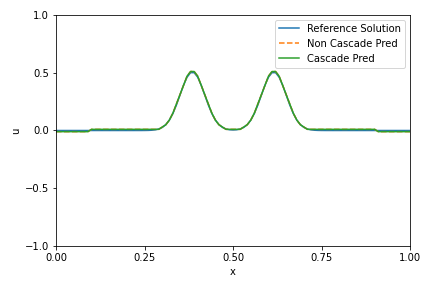}
         \caption{The prediction snapshot at timestep 980 which is in the temporal extrapolation regime, lower triangular encoding.}
     \end{subfigure}
     \hfill
     \begin{subfigure}[t]{0.49\textwidth}
         \centering
         \includegraphics[width=\textwidth]{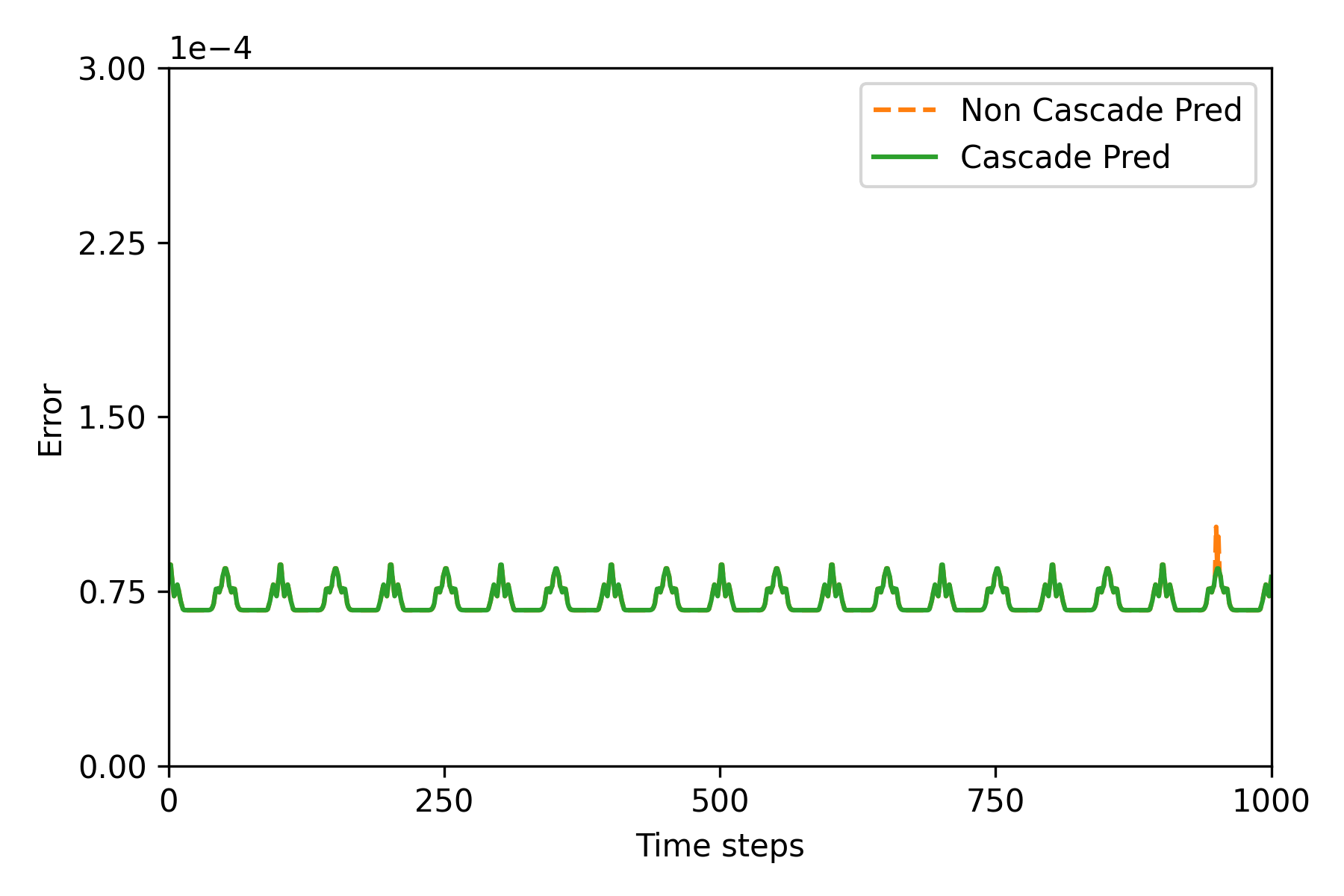}
         \caption{Evolution of error over time-steps. The first 800 timesteps show interpolation and the remaining 200 show extrapolation. Here the error at timestep $n$ is defined as $\frac{1}{N_x}\sum_{j=1}^{N_x}(u_j^n-\Tilde{u}_j^n)^2$, where $u$ is the prediction and $\Tilde{u}$ is the reference solution.}
     \end{subfigure}
     \hfill
     \\
     \begin{subfigure}[t]{0.48\textwidth}
         \centering
         \includegraphics[width=\textwidth]{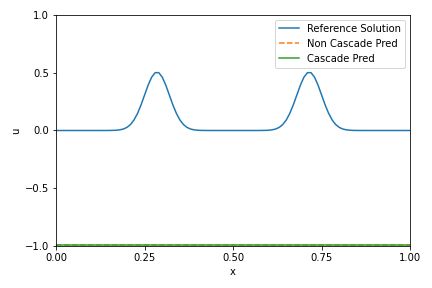}
         \caption{The prediction snapshot at timestep 980 which is in the temporal extrapolation regime, rate encoding.}
     \end{subfigure}
     \hfill
     \begin{subfigure}[t]{0.49\textwidth}
         \centering
         \includegraphics[width=\textwidth]{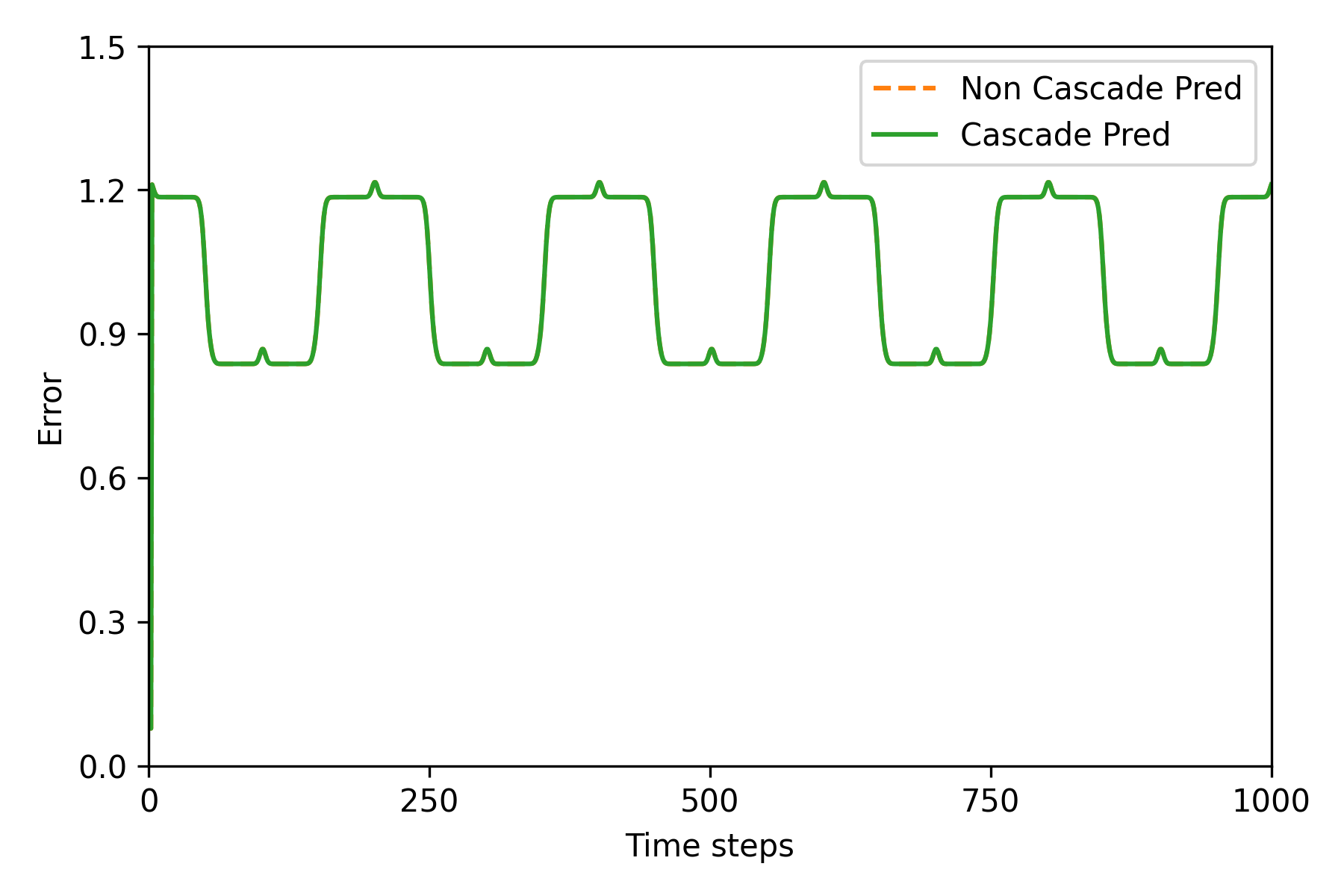}
         \caption{Similar configuration as in (b) but with rate encoding.}
     \end{subfigure}
     \hfill\\
     \begin{subfigure}[t]{0.48\textwidth}
         \centering
         \includegraphics[width=\textwidth]{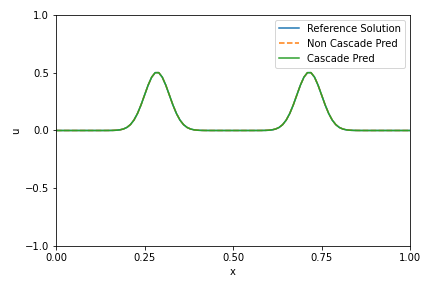}
         \caption{The prediction snapshot at timestep 980 which is in the temporal extrapolation regime, float32 encoding.}
     \end{subfigure}
     \hfill
     \begin{subfigure}[t]{0.49\textwidth}
         \centering
         \includegraphics[width=\textwidth]{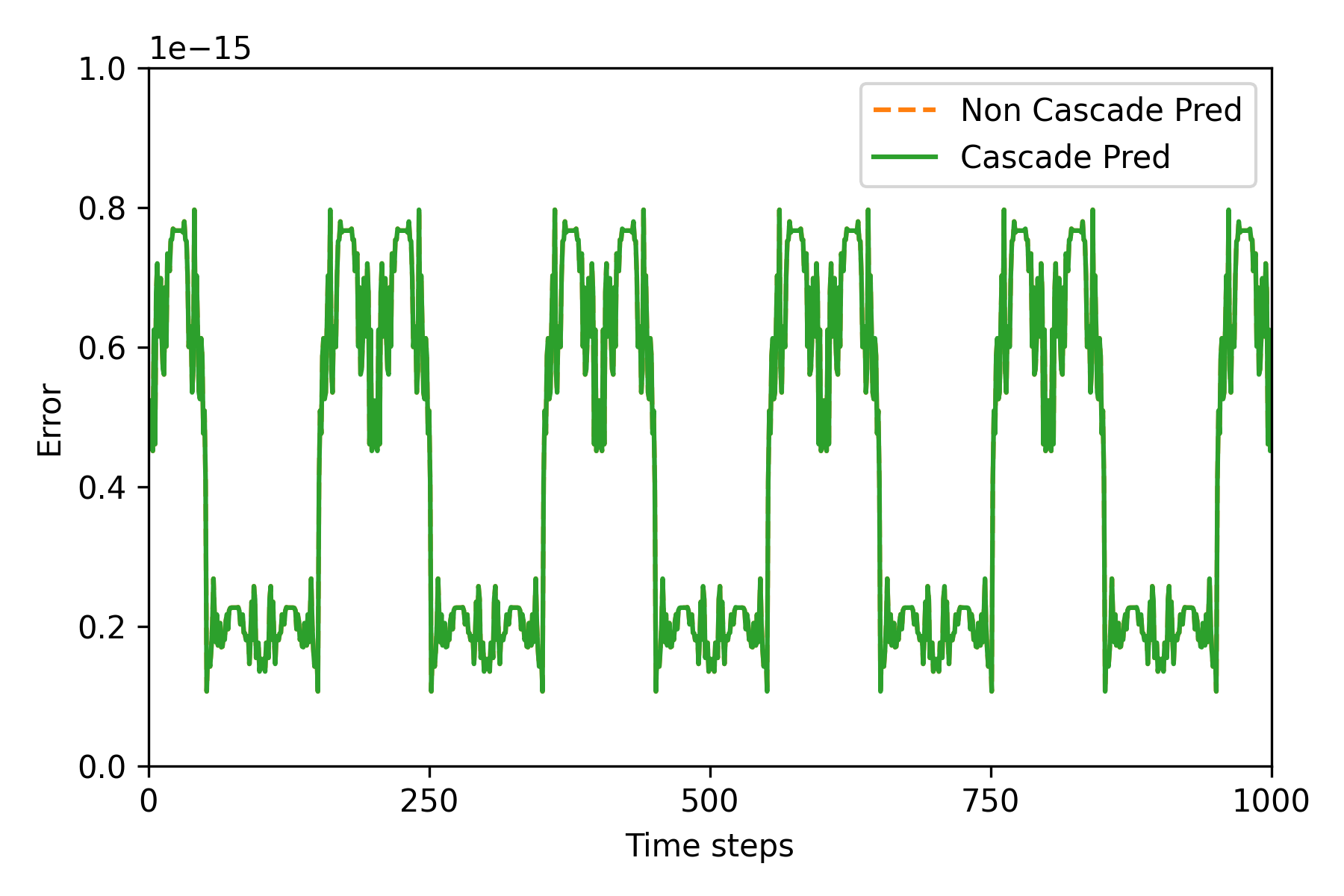}
         \caption{Similar configuration as in (b) but with float32 encoding.}
     \end{subfigure}
     \hfill
     \caption{Wave equation. Results for various encoding schemes.}
     \label{fig:wave}
\end{figure}

We train the network for 5000 epochs, using the first 800 timesteps  for training (interpolation) and the remaining 200 for extrapolation. The results are given in Fig. \ref{fig:wave}. From the results we observe that SMS works well for the wave propagation problem. For this problem we also investigate the effect of the additional two encoding methods. We observe a similar behavior to the one observed in the Lorenz system results. The rate encoding induces a much higher error and produces a completely wrong solution. However, the float32 encoder is performing significantly well, achieving  machine zero error for the cascade prediction at interpolation and the extrapolation intervals.

\subsubsection{Heat equation}

\begin{figure}[!ht]
     \centering
     \begin{subfigure}[t]{0.48\textwidth}
         \centering
         \includegraphics[width=\textwidth]{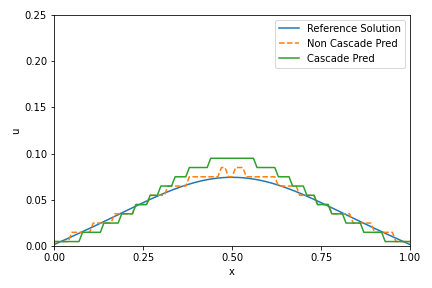}
         \caption{The prediction snapshot at timestep 90 which is in the extrapolation regime.}
     \end{subfigure}
     \hfill
     \begin{subfigure}[t]{0.49\textwidth}
         \centering
         \includegraphics[width=\textwidth]{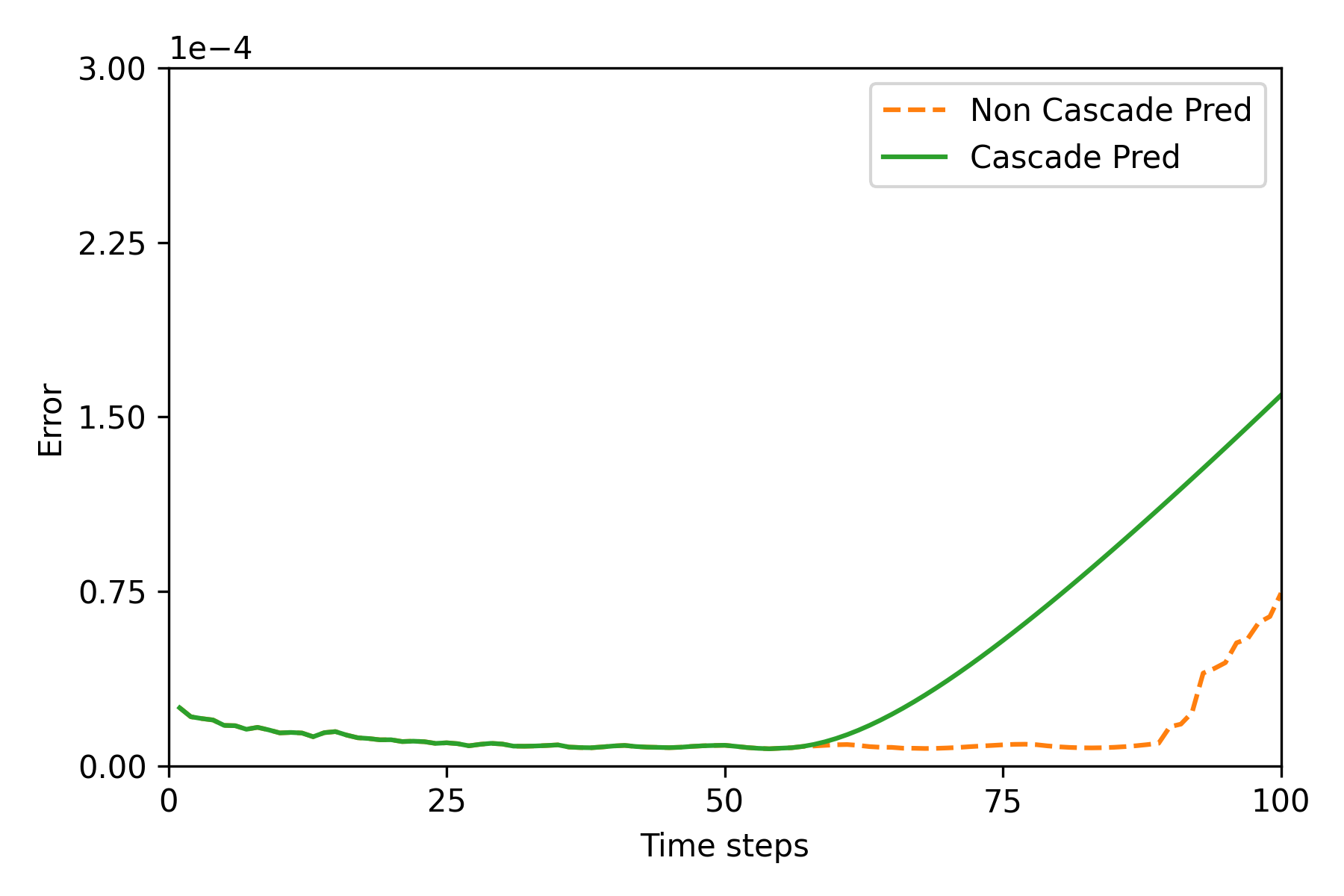}
         \caption{Evolution of error over timesteps. The first 80 time-steps show interpolation and the remaining 20 show extrapolation. Here the error at timestep $n$ is defined as $\frac{1}{N_x}\sum_{j=1}^{N_x}(u_j^n-\Tilde{u}_j^n)^2$, where $u$ is the prediction and $\Tilde{u}$ is the reference solution.}
     \end{subfigure}
     \hfill
     \caption{Heat equation. Encoding method: Lower triangular encoding.}
     \label{fig:heat}
\end{figure}

The one-dimensional heat problem we solve is given by:
\begin{align*}
    \frac{\partial u}{\partial t} = \alpha \frac{\partial^2u}{\partial x^2}&, \qquad t\in [a, b], \quad x\in[x_{min}, x_{max}] \\
    u(0, x) = u_0(x)&, \qquad x\in[x_{min}, x_{max}] \\
\end{align*}
where $u(t, x)$ is the temperature, $x$ is the spatial variable, $t$ is the temporal variable, and $u_0(x)$ is the initial temperature distribution. We use Dirichlet boundary condition set to 0. For this simulation we choose $x_{min} = 0$, $x_{max} = 1$, and the spatial increments are $\Delta x = 0.01$. The diffusion term $\alpha$ is chosen to be 1. The temporal domain is set to $a = 0$ and the temporal increments are chosen to satisfy the condition: $\Delta t = \frac{0.4 * \Delta x ^ 2}{\alpha} = 0.00004$. We integrate for $N_t=1,000$ steps so $b = 0.04$, and use both previous and current time steps to predict the next time step. The number $\#T$ of encoding steps is chosen to be 100 . We sub-sample the data keeping every 10th timestep, so that we have 100 timesteps in total. As an initial condition we choose the function:
\begin{align*}
    f(x) = 
    \begin{cases}
        0&, \quad 0 \leq x < 0.45 \\
        20 (x - 0.45)&, \quad 0.45 \leq x < 0.5 \\
        20 (0.55 - x)&, \quad 0.5 \leq x < 0.55 \\
        0&, \quad 0.55 \leq x \leq 1 \\
    \end{cases}.\label{init_hat}
\end{align*}
We train the network for 5000 epochs and split the data into 80 timesteps used for training and the remaining 20 used for evaluating the extrapolation capabilities. The results are given in Figure \ref{fig:heat}.
We observe that SMS is able to capture the fine changes in the heat diffusion correctly and extrapolate as well (note the small error magnitude).

%% file: conclusion.tex
\section{Discussion and future work}

We have presented a novel method for long time integration using a SNN based marching scheme. We formulated the time integration problem as a data-driven problem and trained a SNN to be an explicit numerical scheme, receiving information about the solution from previous timesteps as input and predicting (inferring) the solution at the next timestep. We conducted several numerical experiments that showed the ability of the proposed method to interpolate and extrapolate, concluding that it is indeed an explicit numerical scheme that can be used for time marching. We measured the error of this spiking marching scheme and showed that the errors stay relatively small, allowing for accurate predictions for long time intervals. To exploit the full potential of the proposed approach, one can implement this on a neuromorphic chip and make the SNN inference very efficient, thus enabling very fine timestepping, leading to highly accurate solutions over time. This can be done in future work when Loihi-2 will be broadly available.

%% file: acknowledgments.tex
\section{Acknowledgements}

This work was supported by the Vannevar Bush Faculty Fellowship award (GEK) from ONR (N00014-22-1-2795).
The work of PS and GEK is supported by the U.S. Department of Energy, Advanced Scientific Computing Research program, under the Scalable, Efficient and Accelerated Causal Reasoning Operators, Graphs and Spikes for Earth and Embedded Systems (SEA-CROGS) project, DE-SC0023191. Pacific Northwest National Laboratory (PNNL) is a multi-program national laboratory operated for the U.S. Department of Energy (DOE) by Battelle Memorial Institute under Contract No. DE-AC05-76RL01830. 